%% file: main.tex
\documentclass[letterpaper, 10 pt, conference]{ieeeconf}  %
\IEEEoverridecommandlockouts
\usepackage[T1]{fontenc}

\title{\LARGE \bf
	Periodic fixed-points and their algebraic characteristics in discrete-time Lur'e feedback systems \thanks{The project was supported by the Israel Science Foundation (grant no. 2406/22) and the Bernard M. Gordon Center for Systems Engineering at the Technion -- IIT, while the second author was also a Jane and Larry Sherman Fellow. All the authors acknowledge the EuroTech Alliance Research Project Initiatives 2025 \& 2026 for funding this research.}
}

\author{Kang Tong, Christian Grussler and Michelle S. Chong   \thanks{K. Tong and C. Grussler are with the Stephen B. Klein Faculty of Aerospace Engineering, Technion -- Israel Institute of Technology, Haifa, Israel. Emails: {\tt\small kang.tong@campus.technion.ac.il, cgrussler@technion.ac.il} }
	\thanks{M. Chong is with the Department of Mechanical Engineering, Eindhoven University of Technology, Eindhoven, Netherlands. Email: {\tt\small m.s.t.chong@tue.nl}}
	\thanks{This paper has been accepted as a regular paper at the 65th IEEE Conference on Decision and Control (CDC 2026).}
}

\include{usr_def}

\begin{document}
	
	\maketitle
	\thispagestyle{empty}
	\pagestyle{empty}

	\begin{abstract}

		We study the problem of identifying nontrivial, i.e., nonzero, periodic fixed-points in discrete-time Lur'e feedback systems. Using the circulant matrix constructed from the transfer function of the linear subsystem, whether stable or unstable, we introduce an algebraic framework that allows us to determine when such fixed-points exist. This framework yields a sector bound defined by two vectors, whose slopes correspond to the maximum and minimum positive singular values of the circulant matrix. 
		Assuming that the nonlinear feedback function is memoryless, we show that a necessary condition for the existence of nontrivial $P$-periodic fixed-points is that the intersection of the continuous completion of the nonlinear feedback function with that sector bound contains at least one point other than the origin.
		Our characterization provides a unified condition valid for all periods $P$, and further enables us to derive upper bounds on the amplitudes of admissible periodic fixed-points with bounded feedback functions. In particular, for relay feedback systems with passive feedback functions, we derive both upper and lower bounds for the amplitudes of such periodic fixed-points.

	\end{abstract}

	\section{Introduction}
	
\begin{figure}[!ht]
		\centering
		\tikzstyle{neu}=[draw, very thick, align = center,circle]
		\tikzstyle{int}=[draw,minimum width=1cm, minimum height=1cm, very thick, align = center]
		\begin{tikzpicture}[>=latex',circle dotted/.style={dash pattern=on .05mm off 1.2mm,
				line cap=round}]
			\node [coordinate, name=input] {};
			\node [int, right of=input, node distance=2cm] (relay) {$\psi(\cdot)$};
			\node [int, right of=relay, node distance=2cm] (system) {$G(z)$};
			\node [coordinate][right of=system, node distance=2cm] (output) {};
			\node [int, below left of=system, node distance=1.5cm, shift={(0cm, -0.5cm)}] (feedback) {$-1$};
			\draw [-] (feedback) -| node[above] {} node[below] {} (input);
			\draw [->] (input) -- node[above] {$u(t)$} (relay);
			\draw [->] (relay) -- node[above] {} (system);
			\draw [-] (system) -- node[above] [name=y] {$y(t)$} (output);
			\draw [->] (output) |- node[above] {} node[below] {} (feedback);
		\end{tikzpicture}
		
		\caption{Feedback system in discrete-time consisting of a linear time-invariant system $G(z)$, $z\in\mathbb{C}$ and a memoryless feedback function $\psi(\cdot)$. \label{fig:lure_feedback}}
	\end{figure}
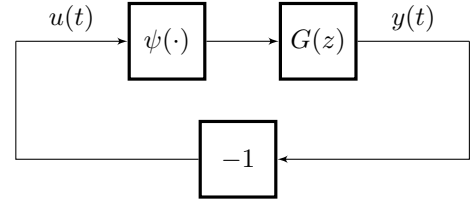

	\emph{Lur'e feedback systems}, as illustrated in \cref{fig:lure_feedback}, which consist of a linear time-invariant (LTI) system in discrete time and a \emph{memoryless} (or \emph{static}) feedback function, permeate many domains, from relay-based PID auto-tuning \cite{aastrom2004revisiting}, power systems \cite{vu_framework_2017}, to neuronal dynamics\cite{iwasaki_lure_2002, gonze_goodwin_2021}.
	In this paper, equilibria and limit cycles in Lur'e feedback systems are characterized as \emph{periodic fixed-points} of an operator derived from the system's loop gain.
	Driven by the practical utility for absolute stability \cite{kim_observer-based_2011}, the prevalence of multiple equilibria \cite{pisarchik_control_2014} and the complexity of oscillatory phenomena \cite{aastrom2004revisiting, gonze_goodwin_2021, miranda-villatoro_analysis_2018}, this study investigates periodic fixed-points.

	Constant periodic fixed-points, i.e., equilibria, in Lur'e feedback systems can be identified in their neighborhoods 
	by Lyapunov stability analysis tools \cite[Chap. 4]{khalil2002nonlinear} and can also indicate global asymptotic stability under appropriate conditions
	\cite{montana_global_2025}.
	Another approach is characterizing the class of nonlinear feedback functions, which ensures that the equilibrium point is globally asymptotically stable in Lur'e systems, known as \emph{the absolute stability problem}. Classical solutions include the small-gain theorem \cite{zames1966input, sontag2002small}, the circle criterion \cite{jayawardhana2011circle}, and the Popov criterion \cite{popov1961absolute, richardson2023strengthened}. Note that the aforementioned conditions were developed for continuous-time systems, and discrete-time counterparts have also been established \cite{jiang2004nonlinear, ahmad_lmi-based_2013, carrasco2019convex}, but they cannot be applied directly to discrete-time (DT) Lur'e feedback systems with multiple constant periodic fixed-points. Additionally, for Lur'e feedback systems with monotone feedback functions, absolute stability can be certified by a Zames–Falb multiplier \cite{turner2021zames, carrasco_zamesfalb_2016, carrasco2019convex}.
	Otherwise, self-oscillations can be identified by constructing destabilizing nonlinearities that fulfill a given sector bound \cite{seiler_construction_2021, kharitenko2024exactness}.

	Nonconstant periodic fixed-points, i.e., limit cycles, have been studied by frequency-domain techniques such as describing function analysis \cite[Chap 7.2]{khalil2002nonlinear} and extensions of the Nyquist criterion for $p$-dominant Lur'e feedback system \cite{miranda-villatoro_analysis_2018}.
	On the other hand, there are also state-space approaches such as Poincaré maps \cite{gonccalves2001global} and the Poincaré--Bendixson theorem for planar systems \cite[Lemma 2.1]{khalil2002nonlinear}.
	However, the generalizations \cite{mallet1990poincare, sanchez_cones_2009, weiss_generalization_2021} require that the feedback functions be monotone and differentiable.

	For the DT Lur'e feedback systems that are widely used in digital control, the literature described thus far has been developed for either equilibria \cite{jiang2004nonlinear, ahmad_lmi-based_2013, kharitenko2024exactness} or limit cycles \cite{mallet1990poincare, sanchez_cones_2009}. This motivates us to develop an alternative viewpoint for providing a structural characterization of system
	behaviors for DT Lur'e feedback systems, which characterizes periodic fixed-points by their periods. We introduce an algebraic framework using characteristics of the circulant matrix \cite[Chap. 3]{gray2005toeplitz}, \cite[Chap. 5]{horn2012matrix} derived from both stable and unstable discrete-time linear time-invariant (DT-LTI) systems to rigorously describe the norm-based bounds for nontrivial $P$-periodic fixed-points, which do not require continuity or differentiability of the feedback functions like in \cite{mallet1990poincare, sanchez_cones_2009, weiss_generalization_2021}.
	For $P=1$, our analysis characterizes nonzero equilibria; for $ P > 1$, it characterizes various limit cycles depending on their periods.
	Furthermore, for all $P \geq 1$, our analysis yields structural conditions that complement existing tools for discrete-time absolute stability analysis, such as small-gain theorems~\cite{jiang2004nonlinear} and stability criteria based on linear matrix inequalities
	(LMIs)~\cite{ahmad_lmi-based_2013}.

	Concretely, we begin our investigations by studying the frequency response of periodic sequences for DT-LTI transfer functions \cite[Chap. 4.4.3]{proakis1996dsp}, and then derive the circulant matrix-based equations for input--output periodic sequences over one period in the time domain, which also apply to unstable DT LTI systems.
	Based on algebraic derivations, we derive a norm-based bound for the feedback function using the maximum and minimum positive singular values of the circulant matrix when the system admits a nontrivial periodic fixed-point.
	Subsequently, we demonstrate the tightest sector bound for the feedback function via the norm-based bound: the system admits nontrivial fixed-points only if continuous completion of its feedback function intersects this sector at a point other than the origin.
	Furthermore, when the feedback function is bounded, we derive an upper bound for the amplitude of nontrivial $P$-periodic fixed-points. Specifically, both upper and lower bounds are provided when the feedback function is a relay with a dead zone. These conclusions are illustrated by numerical examples.

	The remainder of the paper is organized as follows. After some notations and preliminaries in \cref{sec:prelim}, we state our problem in \cref{sec:prob} and provide the algebraic characteristics of the proposed circulant matrix in \cref{sec:circ_mat_alge}. Subsequently, in \cref{sec:main}, we derive our main results on the characterization of nontrivial periodic fixed-points with the feedback function passing through the origin, and illustrate them with examples.
	Finally, \cref{sec:conclusion} draws conclusions, with all proofs deferred to the Appendix.

	\section{Preliminaries} \label{sec:prelim}
	In the following, we introduce several notations and concepts needed for our subsequent derivations and discussions. 
	\subsection{Notations}
	\subsubsection{Sets} We write $\mathbb{R}$ for the set of reals, $\mathbb{C}$ for the set of complex numbers and $\mathbb{Z}$ for the set of integers with $\mathbb{R}_{\geq0}$ ($\mathbb{R}_{>0}$) and $\mathbb{Z}_{\geq0}$ ($\mathbb{Z}_{>0}$) denote the respective subsets of nonnegative (positive) elements. For $k, l \in \mathbb{Z}$ with $k \leq l$, we write $(k:l) := \{ 
	k, k+1, \cdots, l \}$.
	\subsubsection{Sequences} For a sequence $x : \mathbb{Z} \to \mathbb{R}$, with $\sum_{i \in \mathbb{Z}} |x(i)| < \infty$ or $\sup_{i \in \mathbb{Z}} |x(i)| < \infty$, we write $x \in \ell_{1}$ or $x \in \ell_{\infty}$, respectively. 
	For a slice of the sequence $x$, we define the vector $x(k:l):= \begin{bmatrix}
		x(k) & x(k+1) & \dots & x(l)
	\end{bmatrix}^\top$ for $k \leq l$. 
	If there exists a $P \in \mathbb{Z}_{>0}$ such that $x(i) = x(i+P)$ for all $i$, then $x$ is called \emph{$P$-periodic.} The set of all bounded $P$-periodic sequences is denoted by $\ell_\infty(P)$.
	The $P$-truncation vector of $x \in \ell_\infty(P)$ over a single period is denoted by
	\begin{align*}
		x^P := \begin{bmatrix}
			x(0) \\ \vdots \\ x(P-1)
		\end{bmatrix}.
	\end{align*}
	
	\subsubsection{Vectors and matrices} 
	
	For matrices $M = (m_{ij}) \in \mathbb{C}^{n \times m}$, the rank of $M$ is denoted by $\rk(M)$, the image space of $M$ is denoted by $\im(M)$ and the kernel space of $M$ is denoted by $\ker(M)$.
	If $\rk(M) = r \leq \min\{n,m\}$, then $\sigma(M) = \{\sigma_1(M), \sigma_2(M), \cdots, \sigma_r(M)\}$ denotes its \emph{positive singular values}, where the singular values are sorted in descending order, i.e., $\sigma_1(M)$ is the maximum singular value. 
	Let $M^\top$ and $M^*$ denote the transpose and conjugate transpose of $M$.
	Moreover, $M$ is called a \emph{normal matrix} if $M^* M = M M^*$.
	The Moore-Penrose inverse (pseudoinverse) of matrix $M \in \mathbb{C}^{n \times m}$ is denoted by $M^{\dag} \in \mathbb{C}^{m \times n}$. For the vector subspaces $U, V \subseteq \mathbb{R}^n$, $U \perp V$ if and only if $u \perp v$ for all $u\in U$ and $v \in V$.
	
	Let $I_n$ be the identity matrix in $\mathbb{R}^{n \times n}$ and $\boldsymbol{c}_n$ denote the vector with constant entries $c$ in $\mathbb{C}^n$. For a vector $v = [v_1, v_2, \cdots, v_n]^{\top} \in \mathbb{C}^n$, the \emph{circulant matrix} generated from $v$ is denoted by
	\begin{align*}
		H_v := \begin{bmatrix}
			v_1 & v_n & \cdots & v_2 \\
			v_2 & v_1 & \cdots & v_3 \\
			\vdots & \vdots &  & \vdots \\
			v_n & v_{n-1} & \cdots & v_1
		\end{bmatrix} \in \mathbb{C}^{n \times n},
	\end{align*}
	and the \emph{diagonal matrix} generated from $v$ is denoted by
	\begin{align*}
		\diag(v) := \begin{bmatrix}
			v_1 &  &  & \\
			& v_2 &  & \\
			&  & \ddots & \\
			& & & v_n
		\end{bmatrix} \in \mathbb{C}^{n \times n}.
	\end{align*}
	Finally, we denote \emph{the $p$-norm} for vector $v$ with $p\geq 1$ by 
	\begin{align*}
		\| v \|_p := \left( \sum_{i=1}^{n} |v_i|^p \right)^{\frac{1}{p}} \quad \text{and} \quad \| v \|_\infty := \max_i |v_i|.
	\end{align*}
	Subsequently, \emph{the induced matrix $p$-norm} by the vector $p$-norm is defined as
	\begin{align*}
		\| M \|_{p} := \sup_{v \neq 0} \frac{\| M v \|_{p}}{\| v \|_{p}}.
	\end{align*}
	
	\subsubsection{Functions} 
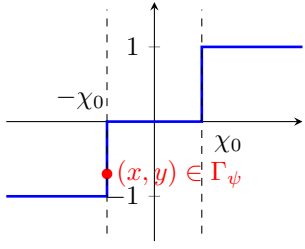
\begin{figure}[!ht]
		\centering
		\begin{tikzpicture}[>=latex']

			\begin{axis}[xtick={0},ytick={-1,1}, width = 5.5  cm, axis lines = middle, 
				ymin = -1.6,    ymax = 1.6,      %
				]
				\addplot[line width = 1 pt,color = blue] table{
					-2.5 -1
					-0.8 -1
					-0.8 0
					0.8 0
					0.8 1
					2.5 1
				};
				\fill[red] (-0.8,-0.7) circle (2pt);
				\node at (-0.8,-0.7) [right, red] {$(x,y) \in \Gamma_{\psi}$};
				
				\draw[dashed] (-0.8,-1.5) -- (-0.8,1.5);
				\draw[dashed] (0.8,-1.5) -- (0.8,1.5);
				\node at (-1.25,0.3){$-\chi_0$};
				\node at (1.25,-0.3){$\chi_0$};
			\end{axis}
		\end{tikzpicture}
		\caption{The continuous completion $\Gamma_\psi$ (the blue curve) when $\psi(\cdot) = \rel_{\chi_0}(\cdot)$ with $\chi_0 \geq 0$.}
		\label{fig:gamma_rel}
	\end{figure}

The \emph{indicator function} of set $\mathcal{A} \subset \mathbb{R}$ is denoted by $\mathbb{1}_{\mathcal{A}}(t)=1$ for $t\in \mathcal{A}$ and $\mathbb{1}_{\mathcal{A}}(t)=0$ for $t \notin \mathcal{A}$.
	In terms of the indicator function $\mathbb{1}_{\mathcal{A}}$, the \emph{unit pulse function} is denoted by $\delta(t):=\mathbb{1}_{\{0\}} (t)$ and the \emph{relay function with symmetric dead zone of width $2\chi_0 \geq 0$} is denoted by
	\begin{align*}
		\rel_{\chi_0}(u) =  - \mathbb{1}_{\mathbb{R}_{< -\chi_0}}(u) + \mathbb{1}_{\mathbb{R}_{> \chi_0}}(u) , \quad u \in \mathbb{R}.
	\end{align*}
	A function $\psi: \mathbb{R} \to \mathbb{R}$ is called \emph{passive} if and only if $\psi(x) x \geq 0$ for all $x \in \mathbb{R}$; its \emph{continuous completion} (a curve) $\Gamma_\psi$ is defined as the smallest closed and connected subset of $\mathbb{R}^2$ that contains the graph of $\psi$. For example, the continuous completion of $\rel_{\chi_0}(\cdot)$ is illustrated in \cref{fig:gamma_rel}.

	\subsection{Linear Time-Invariant Systems} 
	
	We consider finite-dimensional {causal linear discrete-time invariant systems} with (scalar) inputs $u \in \ell_\infty$ and (scalar) outputs $y \in \ell_\infty$. The \emph{impulse response} $g(t) \not \equiv 0$ is the output corresponding to the input $u(t) = \delta(t)$ and the \emph{transfer function} of the system is given by
	\begin{equation} \label{eq:def_trans_fun}
		G(z) = \sum_{t=0}^\infty g(t)z^{-t}, \quad z \in \mathbb{C}.
	\end{equation}
	Additionally, for rational transfer functions, it holds that
	\begin{equation*}
		G(z) = \frac{\beta_m z^m+\beta_{m-1}z^{m-1}+\dots+\beta_1 z + \beta_0}{z^n+\alpha_{n-1}z^{n-1}+\dots+\alpha_1 z + \alpha_0}
	\end{equation*}
	with $n \geq m$ and $\beta_m \neq 0$.

	The \emph{convolution operator} $\mathcal{C}_g: \ell_\infty \to \ell_\infty$ with respect to $G$ is defined by 
	\begin{align} \label{eq:conv_def}
		\mathcal{C}_g(u)(t) := (g \ast u) (t) := \sum_{\tau = -\infty}^{\infty} g(\tau) u(t-\tau), \ t \in \mathbb{Z}.
	\end{align}
	\emph{The discrete Fourier transformation} (DFT) of vector $v \in \mathbb{R}^n$ in frequency domain is defined by $\hat{v} = F_n v$, where
	\begin{align} \label{eq:DFT_mat}
		F_n := \begin{bmatrix}
			1 & 1  & \cdots & 1 \\
			1 & e^{-\frac{2\pi j}{n}}  & \cdots & e^{-\frac{2(n-1)\pi j}{n}} \\
			\vdots & \vdots &  & \vdots \\
			1 & e^{-\frac{2(n-1)\pi j}{n}}  & \cdots & e^{-\frac{2(n-1)^2\pi j}{n}}
		\end{bmatrix}
	\end{align}
	is \emph{the discrete Fourier transformation matrix} of order $n$.
	Subsequently, the \emph{inverse discrete Fourier transformation} (IDFT) is denoted as $v = F_n^{-1} \hat{v} = \frac{1}{n} F_n^{*} \hat{v}$. Those products can be readily calculated by \emph{the fast Fourier transformation} (FFT). 
	
	For $u\in \ell_{\infty}(P)$, it then holds that $y \in \ell_{\infty}(P)$ satisfies $ (F_P y^P) = F_P [\mathcal{C}_g(u)(0), \cdots, \mathcal{C}_g(u)(P-1)]^\top = \Lambda_G^P (F_P u^P)$ by \cite[Chap. 4.4.3]{proakis1996dsp}, where
	\begin{align} \label{eq:Lamda_G^P}
		\Lambda_G^P = \diag([G(1), G(e^{-\frac{2\pi j}{P}}), \cdots, G(e^{-\frac{2(P-1)\pi j}{P}})]).
	\end{align}
	Or equivalently, in the time domain, it holds that
	\begin{align} \label{eq:circul_mat_eq}
		y^P = H_{\overline{g}^P} u^P := F_P^{-1} \Lambda_G^P F_P u^P,
	\end{align}
	where $F_P^{-1} \Lambda_G^P F_P$ is a circulant matrix \cite[Theorem 3.1]{gray2005toeplitz}, which is denoted as $H_{\overline{g}^P}$ with 
	\begin{align} \label{eq:cal_g_P}
		\overline{g}^P := F_P^{-1} \Lambda_G^P \boldsymbol{1}_P.
	\end{align}
	Specifically, for any \emph{bounded-input bounded-output} (BIBO) stable $G(z)$, $\overline{g}^P$ can be computed by the \emph{$P$-periodic summation} of the impulse response of $G(z)$, i.e., $\overline{g}^P_t := \sum_{i=-\infty}^{\infty} g(t + iP + 1)$   $t \in (1:P)$.

	\section{Problem Statement} \label{sec:prob}
	In this work, we study the existence of nontrivial periodic fixed-points in DT Lur'e-feedback systems (see~\cref{fig:lure_feedback}), which consist of a causal linear time-invariant system $G(z)$ with a memoryless nonlinearity $\psi(\cdot)$. Note that the continuity and differentiability of $\psi(\cdot)$ are not required in our analysis.
	
	Using the operator $\mathcal{C}_g$ defined in \cref{eq:conv_def}, this feedback interconnection reads as
	\begin{align} \label{eq:Lure_sys_seq}
		u(t) = -\mathcal{C}_g (\psi(u))(t), \quad t \in \mathbb{Z}.
	\end{align}
	
	We are now ready to precisely define nontrivial periodic fixed-points for \cref{eq:Lure_sys_seq}:
	\begin{defn}[Nontrivial $P$-periodic fixed-point] The sequence $u \in \ell_\infty(P) \setminus \{0\}$ is called a nontrivial $P$-periodic fixed-point if and only if $u$ satisfies \cref{eq:Lure_sys_seq}. Specifically, if a nontrivial $P$-periodic fixed-point $u$ is nonconstant, then $u$ is called a \emph{self-oscillation}.
	\end{defn}

	\section{Circulant Matrix and its algebraic characteristics} \label{sec:circ_mat_alge}
	
	Our main tools are the theory of induced matrix norms \cite[Chap. 5]{horn2012matrix} and circulant matrices \cite[Chap. 3]{gray2005toeplitz}, which allow us to characterize $\psi(\cdot)$ via a sector bound from $H_{\overline{g}^P}$ in \cref{eq:circul_mat_eq}.
	In this section, algebraic properties associated with $H_{\overline{g}^P}$ and the transfer function $G(z)$ are reviewed, including eigenvalues and singular values.

	For the LTI system with transfer function $G(z)$, the eigenvalues denoted by $\lambda(H_{\overline{g}^P})$ can be calculated directly from the corresponding transfer function $G(z)$, which satisfies
	\begin{align} \label{eq:lamda_def}
		\begin{bmatrix}
			\lambda_1(H_{\overline{g}^P}) \\ \lambda_2(H_{\overline{g}^P}) \\
			\vdots \\
			\lambda_P(H_{\overline{g}^P})
		\end{bmatrix} := F_P \overline{g}^P = \begin{bmatrix}
			G(1) \\ G(e^{-\frac{2\pi j}{P}}) \\ \vdots \\ G(e^{-\frac{2(P-1)\pi j}{P}}) 
		\end{bmatrix}.
	\end{align}
	
	Moreover, since $H_{\overline{g}^P}$ is a normal matrix, i.e.,
	\begin{align*}
		H_{\overline{g}^P}^* H_{\overline{g}^P} = (F_P^{-1} \Lambda_G^P F_P)^* (F_P^{-1} \Lambda_G^P F_P) = H_{\overline{g}^P} H_{\overline{g}^P}^*,
	\end{align*} 
	all singular values of $H_{\overline{g}^P}$ are equal to the moduli of the corresponding eigenvalues and \emph{the spectral radius} is denoted as $\rho(H_{\overline{g}^P}) := \max_{n\in (1:P)} | \lambda_n (H_{\overline{g}^P})| = \sigma_1(H_{\overline{g}^P})$.
	Based on the definition of the induced matrix norm, the circulant matrix $H_{\overline{g}^P}$ has the following properties:
	\begin{prope} \label{prope:induced_norm_circ_mat}
		Properties of the induced norms of $H_{\overline{g}^P}$ when $p \geq 1$:
		\begin{enumerate}
			\item $\| H_{\overline{g}^P} \|_2 = \rho(H_{\overline{g}^P}) = \sigma_1(H_{\overline{g}^P})$; \label{prop:induce_norm_1}
			\item $\| H_{\overline{g}^P} \|_{\infty} = \sum_{i=1}^{P} \left| \overline{g}^P_i \right|=\| \overline{g}^P \|_1$; \label{prop:induce_norm_2}
			\item $P^{-\frac{1}{2}}\| H_{\overline{g}^P} \|_{\infty} \leq \| H_{\overline{g}^P} \|_2 \leq P^{\frac{1}{2}} \| H_{\overline{g}^P} \|_{\infty}$. \label{prop:induce_norm_3}
			\item $\| H_{\overline{g}^P} \|_{p} \geq \sigma_1(H_{\overline{g}^P})$ for any $p \geq 1$.
			\label{prop:induce_norm_4}
		\end{enumerate}
	\end{prope}
	
	Additionally, when $\rk(H_{\overline{g}^P})=r \leq P$, the extreme positive singular values satisfy:
	\begin{prope} \label{prope:transfer_fun_to_singular_value}
		When $\rk(H_{\overline{g}^P})=r \leq P$, it implies that
		\begin{align} 
			\sigma_1(H_{\overline{g}^P}) & = \max_{n \in (1:P)} \left\{ \left| G(e^{-\frac{2n\pi j} {P}}) \right| \right\}, \label{eq:sigma_1_P} \\
			\sigma_r(H_{\overline{g}^P}) & = \min_{n \in (1:P)} \left\{ \left| G(e^{-\frac{2n\pi j} {P}}) \right|, G(e^{-\frac{2n\pi j} {P}}) \neq 0 \right\}. \label{eq:sigma_r_P}
		\end{align}
	\end{prope}
	Furthermore, the $\sigma_{\max}$ and the $\sigma_{\min}$ are the extrema of the module of the transfer function $G(z)$ on the unit circle ($z = e^{-j \omega}$ with $\omega \in [0, 2\pi)$), obtained by letting period $P$ tends to infinity:
	\begin{align}
		\sigma_{\max}(G) := \lim_{P\to \infty} \sigma_1(H_{\overline{g}^P}) & = \sup_{\omega \in [0,2\pi)} \left| G(e^{-j\omega}) \right|, \label{eq:max_G}\\
		\sigma_{\min}(G) := \lim_{P\to \infty} \sigma_r(H_{\overline{g}^P}) & = \inf_{\omega \in [0,2\pi)} \left| G(e^{-j\omega}) \right|. \label{eq:min_G}
	\end{align}
	
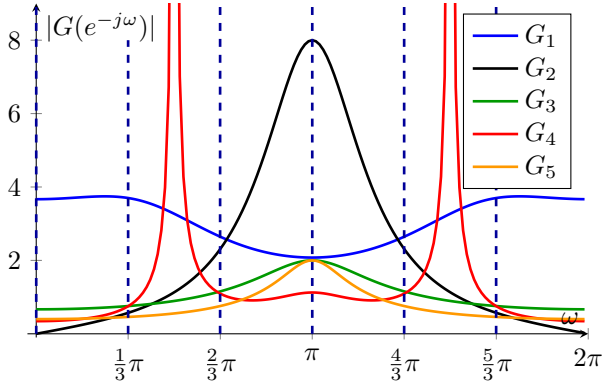
\begin{figure}[!ht]
		\centering
		\begin{tikzpicture}
			\begin{axis}[
				width     = 9cm,
				height    = 6cm,
				xmin = -0.1,    xmax = 2*pi,   
				ymin = -0.1, ymax = 9,
				axis lines = middle,
				legend pos = north east,
				xtick     = {0, 0.333*pi, 0.667*pi, pi, 1.333*pi, 1.667*pi, 2*pi},
				xticklabels = {$0$, $\frac{1}{3}\pi$, $\frac{2}{3}\pi$, $\pi$, $\frac{4}{3}\pi$, $\frac{5}{3}\pi$, $2\pi$},
				xlabel = $\omega$,
				ylabel = $|G(e^{-j\omega})|$,
				]
				\addplot[line width = 1 pt,color = blue] table[col sep=space]  {Nyquist_stable_G_1.txt};
				\addlegendentry{$G_1$};
				\addplot[line width = 1 pt,color = black] table[col sep=space]  {Nyquist_stable_G_2.txt};
				\addlegendentry{$G_2$};
				\addplot[line width = 1 pt,color = black!40!green] table[col sep=space] {Nyquist_stable_G_0.txt};
				\addlegendentry{$G_3$};
				\addplot[line width = 1 pt,color = red] table[col sep=space] {Nyquist_unstable_G_2.txt};
				\addlegendentry{$G_4$};
				\addplot[line width = 1 pt,color = red!40!yellow] table[col sep=space] {Nyquist_very_unstable_G_3.txt};
				\addlegendentry{$G_5$};
				\addplot[dashed, line width = 1 pt, color = black!40!blue] coordinates{(0,0) (0,9)};
				\addplot[dashed, line width = 1 pt, color = black!40!blue] coordinates{(1.046,0) (1.046,8)};
				\addplot[dashed, line width = 1 pt, color = black!40!blue] coordinates{(2.093,0) (2.093,9)};
				\addplot[dashed, line width = 1 pt, color = black!40!blue] coordinates{(3.142,0) (3.142,9)};
				\addplot[dashed, line width = 1 pt, color = black!40!blue] coordinates{(4.188,0) (4.188,9)};
				\addplot[dashed, line width = 1 pt, color = black!40!blue] coordinates{(5.235,0) (5.235,9)};
				
			\end{axis}
		\end{tikzpicture}
		\caption{The magnitudes of five transfer functions $G_1(z)=\frac{3z+0.3}{z^2-0.2z+0.1}$, $G_2(z)=\frac{z-1}{z^2+z+0.25}$, $G_3(z)=\frac{1}{z+0.5}$, $G_4(z)=\frac{z+0.1}{z^3+0.6z^2+z+0.6}$ and $G_5(z)=\frac{1}{z+1.5}$ on the unit circle, i.e., $z = e^{-j\omega}$. The six dashed lines indicate the samples on the unit circle when $P=6$.}
		\label{fig:bode_like}
	\end{figure}

	Note that $\sigma_{\max}(G)$ is known as the $\mathcal{H}_{\infty}$ norm for asymptotically stable LTI systems, which is bounded for the transfer functions that have no poles on the unit circle. As illustrated in \cref{fig:bode_like}, the magnitude of $G_1(z)=\frac{3z+0.3}{z^2-0.2z+0.1}$, $G_2(z)=\frac{z-1}{z^2+z+0.25}$, $G_3(z)=\frac{1}{z+0.5}$, and $G_5(z)=\frac{1}{z+1.5}$ on the unit circle is bounded, while that of $G_4(z)=\frac{z+0.1}{z^3+0.6z^2+z+0.6}$ is unbounded with poles at $\pm j$. 
	Note that the singular values of $H_{\overline{g}^P}$ from the five LTI systems are bounded when $P=6$, which can be read off from the intersections of the plots $|G(e^{-j \omega})|$ with six dashed lines, as shown in \cref{fig:bode_like}.

	For the pseudoinverse of circulant matrix $H_{\overline{g}^P}$, it implies that
	\begin{align*}
		H_{\overline{g}^P}^\dag = (F_P^{-1} {\Lambda_G^P} F_P)^{\dag} = F_P^{-1} {\Lambda_G^P}^{\dag} F_P,
	\end{align*}
	because $F_P^{\dag} = F_P^{-1}$. Note that, since $F_P^{-1} {\Lambda_G^P}^{\dag} F_P$ is also a circulant matrix by \cite[Theorem 3.1]{gray2005toeplitz}, it follows that $H_{\overline{g}^P}^\dag = H_{\tilde{g}^P}$ with
	\begin{align} \label{eq:cal_g_P_pinv}
		\tilde{g}^P = F_P {\Lambda_G^P}^{\dag} \boldsymbol{1}_P.
	\end{align}
	For $\rk(H_{\overline{g}^P}) = r \leq P$, it holds that
	\begin{align} \label{eq:singular_relationship}
		\sigma_1(H_{\overline{g}^P}^{\dag}) = \frac{1}{\sigma_r(H_{\overline{g}^P})}, \quad \sigma_r(H_{\overline{g}^P}^{\dag}) = \frac{1}{\sigma_1(H_{\overline{g}^P})}.
	\end{align}
	
	Finally, using the properties of the image and kernel spaces of $H_{\overline{g}^P}$ and $H_{\overline{g}^P}^\dag$, additional properties based on the norm are derived.
	\begin{prope} \label{prope:kernel_space_pinv}
		For $u^P \in \im(H_{\overline{g}^P})$, $H_{\overline{g}^P}^{\dag} u^P = \boldsymbol{0}_P$ only when $u^P=\boldsymbol{0}_P$.
	\end{prope}
	
	\begin{prope} \label{prope:norm_bound_pinv}
		For $u^P \in \im(H_{\overline{g}^P}) \setminus \{\boldsymbol{0}_P\}$ with $\rk(H_{\overline{g}^P})=r$, it holds that 
		\begin{align} \label{eq:singular_range_pinv}
			\sigma_r(H_{\overline{g}^P}^{\dag}) \| u^P \|_2 \leq \| H_{\overline{g}^P}^{\dag} u^P \|_2 \leq \sigma_1(H_{\overline{g}^P}^{\dag}) \| u^P \|_2,
		\end{align}
		or equivalently, 
		\begin{align}
			\frac{1}{\sigma_1({H_{\overline{g}^P}})} \| u^P \|_2 \leq \| H_{\overline{g}^P}^{\dag} u^P \|_2 \leq \frac{1}{\sigma_r(H_{\overline{g}^P})} \| u^P \|_2.
		\end{align}
	\end{prope}

	\section{Main Results} \label{sec:main}
	
	In this section, we present our main results on nontrivial periodic fixed-points in DT Lur'e feedback systems, where the feedback function is memoryless. We begin by deriving necessary conditions for the existence of nontrivial periodic fixed-points, where the norm of $\psi(u^P)$ is bounded by the norm of $u^P$ multiplied by the extreme positive singular values of $H_{\overline{g}^P}^{\dag}$.
	Based on this, we obtain a sector bound for the passive feedback function $\psi(\cdot)$. For systems with bounded feedback functions, we provide upper bounds for the amplitude of the periodic fixed-point. Furthermore, we derive both lower and upper bounds for the amplitudes in DT relay feedback systems.
	
	\subsection{Necessary conditions for nontrivial periodic fixed-points}

	To analyze $P$-periodic fixed-points in the Lur'e feedback system, we begin by showing that if $u$ satisfies \cref{eq:Lure_sys_seq}, then the periodic fixed-point $u \in \ell_{\infty}(P) \setminus \{\boldsymbol{0}_P\}$ satisfies
	\begin{align} \label{eq:Lure_sys_mat}
		u^P = -H_{\overline{g}^P} \psi(u^P).
	\end{align}
	Note that \cref{eq:Lure_sys_mat} implies the solution $u^P \in \im(H_{\overline{g}^P})$, then for the pseudoinverse form of \cref{eq:Lure_sys_mat}, there exists a vector $u_{\ker} \in \ker(H_{\overline{g}^P})$ such that 
	\begin{align} \label{eq:Lure_sys_mat_pinv}
		\psi(u^P) = -H_{\overline{g}^P}^{\dag} u^P + u_{\ker}.
	\end{align}
	Since ${u^P}^\top u_{\ker} = 0$ by $\im(H_{\overline{g}^P}) \perp \ker(H_{\overline{g}^P})$, it follows that
	\begin{align} \label{eq:Lure_sys_mat_pinv_prod}
		{u^P}^\top \psi(u^P) = {u^P}^\top \left(-H_{\overline{g}^P}^{\dag} u^P \right).
	\end{align}
	Additionally, based on \cref{prope:kernel_space_pinv}, 
	$u^P = \boldsymbol{0}_P$ is a trivial periodic fixed-point if and only if $\psi(0)=0$.
	
	We are now ready to present our first main result: the relationship between the norm of $\psi(u^P)$ and the norm of $u^P$ when $u^P$ is a periodic fixed-point.

	\begin{lem} \label{lem:2_norm_range}
		Let $P \in \mathbb{Z}_{\geq 1}$, $g(t) \not \equiv 0$, $\psi(\cdot)$ be memoryless, and $\rk(H_{\overline{g}^P}) = r$. If \cref{eq:Lure_sys_mat} has a nontrivial $P$-periodic fixed-point, i.e., $u \in \ell_{\infty}(P) \setminus \{\boldsymbol{0}_P\}$, then
		\begin{align} \label{eq:2_norm_range}
			\frac{\| u^P \|_2^2}{\sigma_1({H_{\overline{g}^P}})} \leq \|\psi(u^P)\|_2^2 - \| u_{\ker} \|_2^2 \leq \frac{\| u^P \|_2^2}{\sigma_r({H_{\overline{g}^P}})},
		\end{align} 
		where $u_{\ker} \in \ker(H_{\overline{g}^P})$ satisfying \eqref{eq:Lure_sys_mat_pinv}.
	\end{lem}
	
	It is important to note that the induced $2$-norm yields the smallest gain estimate among all induced $p$-norms when $p \geq 1$, because of $\|H_{\overline{g}^P}\|_p \geq \|H_{\overline{g}^P}\|_2$ from Item 4 in \cref{prope:induced_norm_circ_mat}. Based on this $2$-norm-based bound, we show another necessary condition of $P$-periodic fixed-points by a sector bound.
	
	\begin{thm} \label{thm:necessary_cond_osci}
		Let $P \in \mathbb{Z}_{\geq 1}$, $g(t) \not \equiv 0$, and $\psi(\cdot)$ be memoryless.
		If \cref{eq:Lure_sys_mat} has a nontrivial $P$-periodic fixed-point, then the nonlinear feedback function $\psi(\cdot)$ satisfies 
		\begin{align} \label{eq:2_norm_range_one_side}
			\frac{|x|}{\sigma_1({H_{\overline{g}^P}})} \leq |\psi(x)|,
		\end{align}
		for at least one $x \in \mathbb{R} \setminus \{0\}$.
		
		If $\psi(\cdot)$ is passive additionally, then the continuous completion of $\psi$ satisfies  
		\begin{align} \label{eq:2_norm_range_pos}
			\frac{x^2}{\sigma_1({H_{\overline{g}^P}})} \leq xy \leq \frac{x^2}{\sigma_r({H_{\overline{g}^P}})}
		\end{align}
		for at least one $(x,y) \in \Gamma_\psi$ with $x \in \mathbb{R} \setminus \{0\}$.
	\end{thm}

	In the case of $P=1$ in \cref{lem:2_norm_range} or \cref{thm:necessary_cond_osci}, the nontrivial $1$-periodic fixed-point represents a nontrivial constant sequence, i.e., an equilibrium for the system \cref{eq:Lure_sys_seq}. Then, the \eqref{eq:2_norm_range} reduces to $\|\psi(u^1)\|_2 = \frac{\|u^1\|_2}{|G(1)|}$, or equivalently, $|\psi(u^1)| = \frac{|u^1|}{|G(1)|}$ when $|G(1)| \neq 0$. 
	
	By replacing ${\sigma_1({H_{\overline{g}^P}})}$ and ${\sigma_r({H_{\overline{g}^P}})}$ with $\sigma_{\max}(G)$ and $\sigma_{\min}(G)$ defined in \cref{eq:max_G}-\cref{eq:min_G}, \cref{lem:2_norm_range} and \cref{thm:necessary_cond_osci} provide the common necessary condition for the existence of all nontrivial $P$-periodic fixed-points with $P \geq 1$. For example, \cref{eq:2_norm_range_pos} in \cref{thm:necessary_cond_osci} is represented as 
	\begin{align*}
		\frac{x^2}{\sigma_{\max}(G)} \leq xy \leq \frac{x^2}{\sigma_{\min}(G)}.
	\end{align*}
	If $G(z)$ has no poles and no zeros on the unit circle, then $\sigma_{\max}(G) < \infty$ and $\sigma_{\min}(G)>0$ hold. If $\sigma_{\min}(G)=0$, then define $\frac{x^2}{\sigma_{\min}(G)}=\infty$.
	
	Specifically, based on the DT small gain theorem \cite{jiang2004nonlinear}, every solution of the Lur'e feedback system \cref{eq:Lure_sys_seq} with stable $G(z)$ converges to zero when $|\psi(x)| < \frac{|x|}{\sigma_{\max}(G)}$. In this case, \cref{eq:Lure_sys_seq} does not admit any nontrivial periodic fixed-point, which is consistent with \cref{eq:2_norm_range_one_side} in \cref{thm:necessary_cond_osci}. 
	
	In contrast to the phase condition considered in \cite[Theorem 1]{seiler_construction_2021}, which constructs a nonlinearity $\psi$ that fulfills the given sector bound, \cref{thm:necessary_cond_osci} derives the sector bound condition for the non-trivial $P$-periodic fixed-point from the linear part.

	\subsection{The amplitude of periodic fixed-points}
	
	Next, we will show the bound for the amplitude of periodic fixed-points when the feedback function is bounded, i.e., $\psi(\cdot)$ satisfies \cref{ass:psi_cond_bound}. In the case of relay feedback systems, which are important for relay-based PID auto-tuning, both upper and lower bounds are provided for each period.

	\begin{ass}\label{ass:psi_cond_bound}
		The memoryless feedback function $\psi$ has a bounded range defined by $\mathcal{B}(\psi) := \sup_{x}|\psi(x)| < \infty$.
	\end{ass}
	
	\begin{thm} \label{thm:amplitude_bound}
		Let {$P \in \mathbb{Z}_{\geq 1}$ and memoryless $\psi(\cdot)$ satisfy \cref{ass:psi_cond_bound}.} If \eqref{eq:Lure_sys_mat} has a nontrivial $P$-periodic fixed-point, then the amplitude of $u^P$, i.e., $\|u^P\|_{\infty}$, satisfies
		\begin{align} \label{eq:bound_psi}
			\|u^P\|_{\infty} \leq \| \overline{g}^P \|_1 \mathcal{B}(\psi),
		\end{align}
		or approximately, 
		\begin{align*}
			\|u^P\|_{\infty} \leq P^{\frac{1}{2}} \sigma_1(H_{\overline{g}^P}) \mathcal{B}(\psi).
		\end{align*}
	\end{thm}
	
	Note that if $\sigma_1(H_{\overline{g}^P})=\infty$, i.e., $\Lambda_G^P$ in \cref{eq:circul_mat_eq} contains $\infty$, then \eqref{eq:bound_psi} in \cref{thm:amplitude_bound} is defined by $\|u^P\|_{\infty} \leq \infty$.
	Additionally, the second bound for \cref{thm:amplitude_bound} demonstrates that the amplitude grows no faster than $P^{\frac{1}{2}} \sigma_{\max}(G)$ as the period $P$ increases. 
	
	\begin{cor} \label{cor:amplitude_bound_rfs}
		Let {$P \in \mathbb{Z}_{\geq 1}$}. If the relay feedback system {$u^P = -H_{\overline{g}^P} \rel_{\chi_0}(u^P)$ with $\chi_0 \geq 0$ }has a nontrivial $P$-periodic fixed-point, then {its amplitude $\|u^P\|_{\infty}$} satisfies
		\begin{align} \label{eq:bound_for_RFS}
			\| \tilde{g}^P \|_1^{-1} \leq \|u^P\|_{\infty} \leq \| \overline{g}^P \|_1,
		\end{align}
		where $\tilde{g}^P$ is defined in \eqref{eq:cal_g_P_pinv}, 
		or approximately, 
		\begin{align*}
			P^{-\frac{1}{2}} \sigma_r(H_{\overline{g}^P}) \leq \|u^P\|_{\infty} \leq P^{\frac{1}{2}} \sigma_1(H_{\overline{g}^P}).
		\end{align*}
	\end{cor}
	
	Note that if $\sigma_1(H_{\overline{g}^P})=\infty$, i.e., $\Lambda_G^P$ in \cref{eq:circul_mat_eq} contains $\infty$, then \eqref{eq:bound_for_RFS} in \cref{cor:amplitude_bound_rfs} is defined by $0 \leq \|u^P\|_{\infty} \leq \infty$.

	The following example illustrates \cref{cor:amplitude_bound_rfs} by showing the amplitudes of periodic fixed-points and their bounds.
	
	\begin{exm}
		Consider a DT relay feedback system \cref{eq:Lure_sys_seq} where {$\psi(\cdot)= \rel_{\chi_0}(\cdot)$}. We choose five different transfer functions: $G_1(z) = \frac{3z+0.3}{z^4-0.2z^3+0.1z^2}$, $G_2(z)=\frac{z-1}{z^4+z^3+0.25z^2}$, $G_3(z) = \frac{-1}{z+0.5}$, $G_4(z) = \frac{-z-0.1}{z^3+0.6z^2+z+0.6}$, and $G_5(z) = \frac{-1}{z+1.5}$, whose magnitudes on the unit circle {are the same as those in} \cref{fig:bode_like}. {The configured dead zones $\chi_0$ are listed in \cref{tab:osci_bounds}, together with their amplitudes and the upper and lower bounds calculated from \cref{cor:amplitude_bound_rfs}.
			Subsequently, \cref{fig:five_osci} displays their nontrivial periodic fixed-points, whose periods are 6, 4, 1, 10, and 2. 
			All amplitudes remain within our proposed bounds.}

		\begin{figure}[!ht]
			\centering
			\begin{tikzpicture}[scale=0.9]
				\begin{axis}[height=8 cm,
					width= 9cm,
					xlabel style={right},
					xmax=10.5,xmin = -0.5,
					ymax = 7,ymin=-5,
					xlabel = {$t$},
					ylabel = {$u(t)$},
					axis lines = middle,
					]							
					\addplot+[ycomb,blue,thick, mark=o]
					coordinates{
						(0,-2.3361)
						(1,-4.1570)
						(2,-3.8977)
						(3,2.3361)
						(4,4.1570)
						(5,3.8977) 
						(6,-2.3361) 
						(7,-4.1570) 
						(8,-3.8977)
						(9,2.3361)
						(10,4.1570)
					};
					\addlegendentry{$u_1(t)$}; \label{u1}
					\addplot+[ycomb,black,thick, mark=square]
					coordinates{
						(0,-1.28)
						(1,-0.96)
						(2,1.28)
						(3,0.96)
						(4,-1.28)
						(5,-0.96) 
						(6,1.28) 
						(7,0.96) 
						(8,-1.28)
						(9,-0.96)
						(10,1.28)
					}; 
					\addlegendentry{$u_2(t)$};\label{u2}
					
					\addplot+[ycomb,black!40!green,thick, mark=star]
					coordinates{
						(0,-0.6667)
						(1,-0.6667)
						(2,-0.6667)
						(3,-0.6667)
						(4,-0.6667)
						(5,-0.6667) 
						(6,-0.6667) 
						(7,-0.6667) 
						(8,-0.6667)
						(9,-0.6667) 
						(10,-0.6667)
					}; 
					\addlegendentry{$u_3(t)$};\label{u3}
					
					\addplot+[ycomb,red,thick, mark=triangle]
					coordinates{
						(0,-0.3676)
						(1,-0.4554)
						(2,-0.1507)
						(3,-0.5855)
						(4,-0.0726)
						(5,0.3676)
						(6,0.4554)
						(7,0.1507)
						(8,0.5855)
						(9,0.0726)
						(10,-0.3676)
					}; 
					\addlegendentry{$u_4(t)$};\label{u4}
					
					\addplot+[ycomb,red!40!yellow,thick, mark=diamond]
					coordinates{
						(0,-2)
						(1,2)
						(2,-2)
						(3,2)
						(4,-2)
						(5,2) 
						(6,-2) 
						(7,2) 
						(8,-2)
						(9,2) 
						(10,-2)
					}; 
					\addlegendentry{$u_5(t)$};\label{u5}
					
				\end{axis}
			\end{tikzpicture}
			\caption{Nontrivial periodic fixed-points in five different relay feedback systems with transfer functions \ref{u1} $G_1(z) = \frac{3z+0.3}{z^4-0.2z^3+0.1z^2}$, \ref{u2} $G_2(z)=\frac{z-1}{z^4+z^3+0.25z^2}$, \ref{u3} $G_3(z)=\frac{-1}{z+0.5}$, \ref{u4} $G_4(z) = \frac{-z-0.1}{z^3+0.6z^2+z+0.6}$, and \ref{u5} $G_5(z) = \frac{-1}{z+1.5}$ in Example 1.}
			\label{fig:five_osci}
		\end{figure}
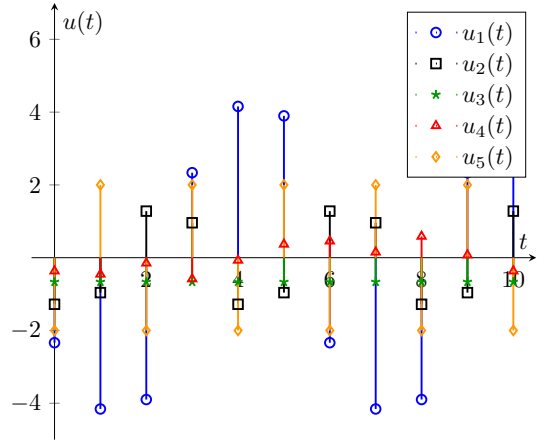
		
		\begin{table}[htbp]
			\centering
			\caption{The periods, amplitudes, and upper and lower bounds for nontrivial periodic fixed-points in Example 1.}
			\label{tab:osci_bounds}
			
			\begin{tabular}{l *{5}{S[table-format=2.3]}}
				\toprule
				\textbf{} & {$G_1(z)$} & {$G_2(z)$} & {$G_3(z)$} & {$G_4(z)$} & {$G_5(z)$} \\
				\midrule
				{Dead zone $\chi_0$} & 0 & 0 & 0.1 & 0.15 & 1.5 \\
				Period  & 6 & 4 & 1 & 10  & 2 \\
				Upper bound   & 24.94 & 32.00 & 1.50 & 31.30 & 4.00 \\
				Lower bound   & 0.42 & 0.48 & 0.67 & 0.37 & 0.96 \\
				Observed value & 4.16 & 1.28 & 0.67 & 0.59 & 2.00 \\
				\bottomrule
			\end{tabular}    
		\end{table}
		
	\end{exm}

	{Note that $G_4(z)$ yields a finite bound for $10$-periodic fixed-points even when $\sigma_{max}(G_4) = \infty$ as shown in \cref{fig:bode_like}.
		Although $G_5(z)$ is an unstable system, the resulting self-oscillation still fulfills the amplitude bound \cref{eq:bound_for_RFS} in \cref{cor:amplitude_bound_rfs}.
	}

	\section{Conclusion} \label{sec:conclusion}
	
	This paper has investigated the existence of nontrivial {periodic} fixed-points and their amplitudes in DT Lur'e feedback systems with a {memoryless} feedback function. By exploiting algebraic characteristics of the circulant matrix, we derived that the norm of the {nontrivial} periodic fixed-point is bounded by the maximum and minimum positive singular values of the constructed circulant matrix, which can be directly derived from the transfer function.
	Based on these insights, we derived necessary conditions for the existence of nontrivial periodic fixed-points and provided explicit amplitude bounds for systems with bounded nonlinearities, including relay feedback.
	We anticipate that this discrete-time analysis will provide fundamental insights into the characteristics of {nontrivial} periodic fixed-points, {identifying the structural characterization} of equilibria and limit cycles.
	Additionally, our analysis also applies to unstable linear systems. The stability of {our} identified {nontrivial} periodic fixed-points is not addressed in this paper.
	Future research will aim to {establish necessary and sufficient conditions for their existence, thereby providing a structural characterization of equilibria and limit cycles.}

	\appendix

	\subsection{Proof of \cref{prope:induced_norm_circ_mat}}
	
	From \cite[Example 5.6.6.]{horn2012matrix}, it has that $\| H_{\overline{g}^P} \|_2 = \sigma_1(H_{\overline{g}^P})$. Since $H_{\overline{g}^P}$ is a normal matrix, it implies $\sigma_1(H_{\overline{g}^P}) = \rho(H_{\overline{g}^P})$, i.e.,
	\cref{prop:induce_norm_1} in \cref{prope:induced_norm_circ_mat} holds.
	
	From \cite[Example 5.6.5.]{horn2012matrix}, it has that 
	\begin{align*}
		\| H_{\overline{g}^P} \|_{\infty} = \max_{i \in (1:P)} \sum_{j=1}^{P}|(H_{\overline{g}^P})_{ij}| = \sum_{k=1}^P |{\overline{g}_k^P}| = \| \overline{g}^P \|_1,
	\end{align*}
	where $(H_{\overline{g}^P})_{ij}$ denotes the i-th row and j-th column's element in $H_{\overline{g}^P}$. Thus,  \cref{prop:induce_norm_2} in \cref{prope:induced_norm_circ_mat} holds.
	
	The \cref{prop:induce_norm_3} in \cref{prope:induced_norm_circ_mat} is proved by \cite[5.6.P23]{horn2012matrix}.
	
	From \cite[Thm. 5.6.9.]{horn2012matrix}, it holds that $\| H_{\overline{g}^P} \|_{p} \geq \rho(H_{\overline{g}^P})$ for all $p\geq 1$. Since $\rho(H_{\overline{g}^P}) = \sigma_1(H_{\overline{g}^P})$ by \cref{prop:induce_norm_1}, \cref{prop:induce_norm_4} in \cref{prope:induced_norm_circ_mat} holds.

	\subsection{Proof of \cref{prope:transfer_fun_to_singular_value}}
	
	When $\rk(H_{\overline{g}^P})=r \leq P$, $H_{\overline{g}^P}$ has $r$ non-zero eigenvalues. Since $H_{\overline{g}^P}$ is a normal matrix and its eigenvalues are denoted by \cref{eq:lamda_def}, it implies that
	\begin{align*}
		\sigma_1(H_{\overline{g}^P}) & = \max_{n \in (1:P)} \left\{ \left| \lambda_n(H_{\overline{g}^P}) \right| \right\}=  \max_{n \in (1:P)} \left\{ \left| G(e^{-\frac{2n\pi j} {P}}) \right| \right\}, \\
		\sigma_r(H_{\overline{g}^P}) &  = \min_{n \in (1:P)} \left\{ \left| \lambda_n(H_{\overline{g}^P}) \right|, \lambda_n(H_{\overline{g}^P}) \neq 0 \right\} \\
		& = \min_{n \in (1:P)} \left\{ \left| G(e^{-\frac{2n\pi j} {P}}) \right|, G(e^{-\frac{2n\pi j} {P}}) \neq 0 \right\}.
	\end{align*}
	
	\subsection{Proof of \cref{prope:kernel_space_pinv}}
	
	As $u^P \in \im(H_{\overline{g}^P})$, $H_{\overline{g}^P}^{\dag} u^P = \boldsymbol{0}_P$ only when $u^P \in \ker(H_{\overline{g}^P}^{\dag}) \cap \im(H_{\overline{g}^P})$. 
	By $\ker(H_{\overline{g}^P}^{\dag})\perp \im(H_{\overline{g}^P})$ from \cite[Proposition 3.3]{barata_moorepenrose_2012}, it holds that $\ker(H_{\overline{g}^P}^{\dag}) \cap \im(H_{\overline{g}^P}) = \{\boldsymbol{0}_P\}$. \hfill $\Box$
	
	\subsection{Proof of \cref{prope:norm_bound_pinv}}
	
	Letting $\tilde{x} = \frac{u^P}{\| u^P \|_2}$ with $u^P \in \im(H_{\overline{g}^P}) \setminus \{\boldsymbol{0}_P\}$, it holds that $\tilde{x} \perp \ker(H_{\overline{g}^P}^{\dag})$ and $\| \tilde{x}\|_2 = 1$. By singular value decomposition, it follows that $H_{\overline{g}^P}^\dag = U \Sigma V^*$ with $\Sigma = \diag([\sigma_1(H_{\overline{g}^P}^\dag), \cdots, \sigma_r(H_{\overline{g}^P}^\dag), 0, \cdots, 0])$ and $\sigma_1(H_{\overline{g}^P}^\dag) \geq \sigma_2(H_{\overline{g}^P}^\dag) \geq \cdots \geq \sigma_r(H_{\overline{g}^P}^\dag) >0$. Letting $y = V^* \tilde{x} = [y_1, y_2, \cdots, y_P]^{\top}$, then it implies that $\| y \|_2 = \| \tilde{x}\|_2 = 1$ with $y_{r+1} = \cdots = y_P = 0$ by $\tilde{x} \perp \ker(H_{\overline{g}^P}^{\dag})$.
	
	Hence, the range of $\| H_{\overline{g}^P}^{\dag} \tilde{x} \|_2$ is $[\sigma_r(H_{\overline{g}^P}^{\dag}), \sigma_1(H_{\overline{g}^P}^{\dag})]$ by
	\begin{align*}
		\| H_{\overline{g}^P}^{\dag} \tilde{x} \|_2 & = \| U\Sigma y \|_2 = {\sqrt{(U\Sigma y)^{*} (U\Sigma y)}} = \sqrt{y^* \Sigma^* \Sigma y} \\
		& = {\left( \sigma_1(H_{\overline{g}^P}^{\dag})^2 y_1^2 + \cdots + \sigma_r(H_{\overline{g}^P}^{\dag})^2 y_r^2 \right)^{\frac{1}{2}}}.
	\end{align*}
	
	Since $\sum_{i=1}^r y_i^2 = \| y \|_2^2 = 1$, it implies that
	\begin{align*}
		\sigma_r(H_{\overline{g}^P}^{\dag})^2 \leq \left( \sum_{i=1}^{r} \sigma_i(H_{\overline{g}^P}^{\dag})^2 y_i^2 \right)^{\frac{1}{2}} \leq \sigma_1(H_{\overline{g}^P}^{\dag})^2.
	\end{align*}
	Specifically, $\| H_{\overline{g}^P}^{\dag} \tilde{x} \|_2 = \sigma_1(H_{\overline{g}^P}^{\dag})$ when $|y_1| = 1$ and $y_i = 0$ for $i\in(2:r)$; $\| H_{\overline{g}^P}^{\dag} \tilde{x} \|_2 = \sigma_r(H_{\overline{g}^P}^{\dag})$ when $|y_r| = 1$ and $y_i = 0$ for $i\in(1:r-1)$. Since $\| H_{\overline{g}^P}^{\dag} u^P \|_2 = \| H_{\overline{g}^P}^{\dag}  (\|u^P \|_2\tilde{x}) \|_2 = \|u^P \|_2 \| H_{\overline{g}^P}^{\dag} \tilde{x} \|_2$, the inequalities \cref{eq:singular_range_pinv} hold. \hfill $\Box$
	
	\subsection{Proof of \cref{lem:2_norm_range}}
	{Since $\psi(u^P) = -H_{\overline{g}^P}^{\dag} u^P + u_{\ker}$ with $u^P \in \im(H_{\overline{g}^P}) \setminus \{\boldsymbol{0}_P\}$ and $u_{\ker} \in \ker(H_{\overline{g}^P})$, it implies that
		\begin{align*}
			\|\psi(u^P)\|_2^2 & = \| -H_{\overline{g}^P}^{\dag} u^P + u_{\ker} \|_2^2 \\
			& = {u^P}^\top {H_{\overline{g}^P}^{\dag}}^\top H_{\overline{g}^P}^{\dag} u^P + u_{\ker}^\top u_{\ker} \\ 
			& = \| H_{\overline{g}^P}^{\dag} u^P \|_2^2 + \| u_{\ker} \|_2^2.
		\end{align*}
		By the bound from \cref{prope:norm_bound_pinv}, it holds that
		\begin{align*}
			\frac{\| u^P \|_2^2}{\sigma_1({H_{\overline{g}^P}})}+\| u_{\ker} \|_2^2 \leq \|\psi(u^P)\|_2^2 \leq \frac{\| u^P \|_2^2}{\sigma_r({H_{\overline{g}^P}})}+\| u_{\ker} \|_2^2,
		\end{align*} 
		which is equal to the provided inequality in \cref{eq:2_norm_range}.
		
		If $\rk(H_{\overline{g}^P})=P$ additionally, then $u_{\ker} = \boldsymbol{0}_P$ by $\ker(H_{\overline{g}^P}) = \{ \boldsymbol{0}_P \}$; it follows 
		\begin{align*}
			\frac{\| u^P \|_2}{\sigma_1({H_{\overline{g}^P}})} \leq \|\psi(u^P)\|_2 \leq \frac{\| u^P \|_2}{\sigma_P({H_{\overline{g}^P}})}.
		\end{align*} 
		\hfill $\Box$
	}

	\subsection{Proof of \cref{thm:necessary_cond_osci}}
	We prove it by contradiction. 
	{Assuming $|\psi(x)| < \frac{|x|}{\sigma_1({H_{\overline{g}^P}})}$ for all $x \in \mathbb{R} \setminus \{0\}$. For all $u^P \in \im(H_{\overline{g}^P}) \setminus \{\boldsymbol{0}_P\}$, it follows that
		\begin{align*}
			\|\psi(u^P)\|_2 < \frac{\|u^P\|_2}{\sigma_1({H_{\overline{g}^P}})},
		\end{align*}
		which contradicts \cref{lem:2_norm_range}. Hence, it holds that $\frac{|x|}{\sigma_1({H_{\overline{g}^P}})} \leq |\psi(x)|$ for at least one $x$.
		
		Furthermore, when $\psi(\cdot)$ is passive, it has 
		\begin{align*}
			\frac{|x|}{\sigma_1({H_{\overline{g}^P}})} |x| \leq |\psi(x)||x| \iff \frac{x^2}{\sigma_1({H_{\overline{g}^P}})} \leq \psi(x)x.  
		\end{align*}
		Assuming $x\psi(x) > \frac{x^2}{\sigma_r({H_{\overline{g}^P}})}$ for all $x \in \mathbb{R} \setminus \{0\}$ for any $u^P \in \im(H_{\overline{g}^P}) \setminus \{\boldsymbol{0}_P\}$, \eqref{eq:Lure_sys_mat_pinv_prod} implies that 
		{\small\begin{align*}
				& \quad \| {u^P}^\top \left(-H_{\overline{g}^P}^{\dag} u^P \right) \|_2 = \| {u^P}^\top \psi(u^P) \|_2  = \sum_{i=1}^{P} u^P_i \psi(u^P_i)\\
				& \implies \| {u^P}^\top H_{\overline{g}^P}^{\dag} u^P  \|_2 = \sum_{i=1}^{P} u^P_i \psi(u^P_i) \tag*{by passive $\psi(\cdot)$} \\
				& \implies \sum_{i=1}^{P} u^P_i \psi(u^P_i) \leq \sigma_1(H_{\overline{g}^P}) \|u^P\|_2^2 = \sigma_1(H_{\overline{g}^P})\sum_{i=1}^{P} (u^P_i)^2 \tag*{by \cref{prope:norm_bound_pinv}}, 
		\end{align*}}
		which violates $\psi(u^P_i) x > \frac{(u^P_i)^2}{\sigma_r({H_{\overline{g}^P}})}$ by the provided assumption.
		
		Hence, $\Gamma_{\psi}$ can lie neither entirely in $\{ (x,y): xy > \frac{x^2}{\sigma_r({H_{\overline{g}^P}})}\}$ nor entirely in $\{ (x,y): xy < \frac{x^2}{\sigma_1({H_{\overline{g}^P}})}\}$, it must contain a point $(x,y) \in \Gamma_\psi$ with $x \in \mathbb{R} \setminus \{0\}$ satisfying 
		\begin{align*} 
			\frac{x^2}{\sigma_1({H_{\overline{g}^P}})} \leq xy \leq \frac{x^2}{\sigma_r({H_{\overline{g}^P}})}.
	\end{align*}}
	\hfill $\Box$

	\subsection{Proof of \cref{thm:amplitude_bound}}
	Since every nontrivial $P$-periodic fixed-point satisfies $u^P = H_{\overline{g}^P} \psi(u^P)$, it implies that
	\begin{align*}
		\|u^P\|_{\infty} & = \| -H_{\overline{g}^P} \psi(u^P) \|_{\infty} \leq \| H_{\overline{g}^P} \|_{\infty} \| \psi(u^P) \|_{\infty} \\
		& \leq \| H_{\overline{g}^P} \|_{\infty} \mathcal{B}(\psi) = \| \overline{g}^P \|_{1} \mathcal{B}(\psi),
	\end{align*}
	where the second equality follows from \cref{prop:induce_norm_2} in \cref{prope:induced_norm_circ_mat}.
	
	Since $P^{-\frac{1}{2}}\| H_{\overline{g}^P} \|_{\infty} \leq \| H_{\overline{g}^P} \|_{2}$ by \cref{prop:induce_norm_3} in \cref{prope:induced_norm_circ_mat}, it follows that
	\begin{align*}
		\|u^P\|_{\infty} & \leq \| -H_{\overline{g}^P} \|_{\infty} \mathcal{B}(\psi) \leq P^{\frac{1}{2}} \| H_{\overline{g}^P} \|_{2} \mathcal{B}(\psi) \\
		& = P^{\frac{1}{2}} \sigma_1(H_{\overline{g}^P}) \mathcal{B}(\psi).
	\end{align*}
	\hfill $\Box$

	\subsection{Proof of \cref{cor:amplitude_bound_rfs}}
	For a nontrivial $P$-periodic fixed-point, it implies {$\left\| \rel_{\chi_0}(u^P) \right\|_{\infty} = 1$}. Since $\rel_{\chi_0}(\cdot)$ fulfills \cref{ass:psi_cond_bound}, the upper bound for \cref{thm:amplitude_bound} reduces to $\|u^P\|_{\infty} \leq \| \overline{g}^P \|_1$. For the pseudoinverse system \eqref{eq:Lure_sys_mat_pinv}, it also implies $\|\psi(u^P)\|_{\infty} = \| H_{\overline{g}^P}^{\dag} u^P \|_{\infty} \leq \| H_{\overline{g}^P}^{\dag} \|_{\infty} = \| {H_{\tilde{g}^P}} \|_{\infty}$, where $\tilde{g}^P$ is defined in \cref{eq:cal_g_P_pinv}. Since $\| H_{\tilde{g}^P} \|_{\infty} = \| {\tilde{g}^P} \|_1$ by \cref{prop:induce_norm_2} in \cref{prope:induced_norm_circ_mat}, it follows $\|u^P\|_{\infty} \geq \| {H_{\tilde{g}^P}} \|_{\infty}^{-1} = \| \tilde{g}^P \|_1^{-1}$.
	From $\| H_{\overline{g}^P} \|_2 \leq P^{\frac{1}{2}} \| H_{\overline{g}^P} \|_{\infty}$ by \cref{prop:induce_norm_3} in \cref{prope:induced_norm_circ_mat}, the approximated bound is verified. \hfill $\Box$

	\bibliographystyle{ieeetr}
	\bibliography{ref}

\end{document}

%% file: usr_def.tex
\usepackage{graphicx}
\usepackage{enumerate}
\usepackage{mathtools}
\usepackage{xcolor} %

\usepackage{amsmath} %
\usepackage{amssymb}  %
\usepackage{bm}
\usepackage{fixme}
\fxsetup{status=draft}

\usepackage{bbm}
\usepackage{algorithm}
\usepackage{algorithmic}
\usepackage{pgfplots}
\pgfplotsset{compat=newest}

\usepackage{tikz}
\usepackage{tkz-euclide}

\usepackage{enumerate}

\usepackage{mathrsfs}
\definecolor{ao(english)}{rgb}{0.0, 0.5, 0.0}
\usepackage[colorlinks, citecolor = {ao(english)}, linkcolor = {ao(english)}]{hyperref} 
\usepackage[capitalize,nameinlink]{cleveref}
\usepackage{siunitx}
\usepackage[sort,nocompress]{cite}
\usepackage{subcaption}

\usepackage{booktabs}      %
\usepackage{multirow}      %

\usetikzlibrary{arrows,shapes,trees,calc,positioning,patterns,decorations.pathmorphing,decorations.markings,backgrounds}
\usetikzlibrary{matrix}

\usepgfplotslibrary{groupplots}

\newtheorem{thm}{Theorem}
\crefname{thm}{Theorem}{Theorems}

\crefname{prop}{Proposition}{Propositions}

\newtheorem{lem}{Lemma}
\crefname{lem}{Lemma}{Lemmas}

\newtheorem{prope}{Property}
\crefname{prope}{Property}{Properties}

\newtheorem{cor}{Corollary}
\crefname{cor}{Corollary}{Corollaries}

\crefname{rem}{Remark}{Remark}

\newtheorem{ass}{Assumption}
\crefname{ass}{Assumption}{Assumption}

\crefname{conj}{Conjecture}{Conjectures}

\newtheorem{defn}{Definition}
\crefname{defn}{Definition}{Definitions}

\crefname{prob}{Problem}{Problems}

\crefname{appl}{Application}{Applications}

\newtheorem{exm}{Example}
\crefname{exm}{Example}{Examples}

\crefname{algorithm}{Algorithm}{Algorithms}
\crefname{paper}{Paper}{Papers}
\crefname{figure}{Figure}{Figures}
\crefname{section}{Section}{Sections}

\usepackage{dsfont}
\let\mathbb=\mathds

\newcommand{\rk}{\textnormal{rank}}

\newcommand{\diag}{\textnormal{diag}}

\newcommand{\im}{\textnormal{im}}

\newcommand{\rel}{\textnormal{rel}}

\colorlet{FigColor1}{blue}
\colorlet{FigColor2}{red}
\colorlet{FigColor3}{ao(english)}
\colorlet{FigColor4}{orange}
\pgfplotsset{every axis plot/.append style={line width=1.5pt}}
\definecolor{bluebell}{rgb}{0.74, 0.83, 0.9}
\definecolor{airforceblue}{rgb}{0.36, 0.54, 0.66}

\crefformat{equation}{\textup{#2(#1)#3}}
\crefrangeformat{equation}{\textup{#3(#1)#4--#5(#2)#6}}
\crefmultiformat{equation}{\textup{#2(#1)#3}}{ and \textup{#2(#1)#3}}
{, \textup{#2(#1)#3}}{, and \textup{#2(#1)#3}}
\crefrangemultiformat{equation}{\textup{#3(#1)#4--#5(#2)#6}}%
{ and \textup{#3(#1)#4--#5(#2)#6}}{, \textup{#3(#1)#4--#5(#2)#6}}{, and \textup{#3(#1)#4--#5(#2)#6}}

\Crefformat{equation}{#2Equation~\textup{(#1)}#3}
\Crefrangeformat{equation}{Equations~\textup{#3(#1)#4--#5(#2)#6}}
\Crefmultiformat{equation}{Equations~\textup{#2(#1)#3}}{ and \textup{#2(#1)#3}}
{, \textup{#2(#1)#3}}{, and \textup{#2(#1)#3}}
\Crefrangemultiformat{equation}{Equations~\textup{#3(#1)#4--#5(#2)#6}}%
{ and \textup{#3(#1)#4--#5(#2)#6}}{, \textup{#3(#1)#4--#5(#2)#6}}{, and \textup{#3(#1)#4--#5(#2)#6}}

\crefdefaultlabelformat{#2\textup{#1}#3}